\documentclass[graybox]{svmult}
\usepackage{mathptmx}       
\usepackage{helvet}         
\usepackage{courier}        
\usepackage{makeidx}         
\usepackage{graphicx}        
\usepackage{multicol}        
\usepackage[bottom]{footmisc}
\usepackage{amsmath,amssymb,amsfonts}        
\usepackage{xcolor}
\makeindex             

\newtheorem{prop}{Proposition}

\begin{document}
\title*{Delay differential equations with periodic coefficients: a numerical insight}
\author{Anatoli Ivanov and Sergiy Shelyag}
\institute{Anatoli Ivanov \at Penn State Wilkes-Barre, 44 University Drive, Dallas, PA 18612, USA \\ \email{afi1@psu.edu}
\and Sergiy Shelyag \at College of Science and Engineering, Flinders University at Tonsley, Adelaide, South Australia 5042, Australia \\ \email{sergiy.shelyag@flinders.edu.au}}
%
%
\maketitle

\abstract{
Simple form scalar differential equation with delay and non-linear negative periodic feedback is considered. 
The existence of slowly oscillating periodic solutions with the same period as the feedback coefficient is shown numerically within 
the admissible range for the periods.
}

\section{Introduction}
\label{sec:1}

This paper presents sample data obtained from the analysis of numerical solutions of scalar simple form delay differential equations (DDEs) with non-linear negative periodic feedback. The aim of this preliminary analysis is to study periodic solutions of these equations and their dynamics in general. A general form of the class of DDEs we aim our research efforts at is given by the following: 
\begin{equation}\label{DDE-Per}
    x^\prime(t)=-\mu(t) x(t)+a(t)f(x(t-\tau)),\quad \tau>0,
\end{equation}
where the non-linearity $f(x)$ is a continuous function and the coefficients $a(t)$ and $\mu(t)$ are continuous  $T$-periodic functions. 
This equation is a natural generalization to the periodic case of the well known 
 autonomous equation
\begin{equation}\label{DDE-Aut}
    x^\prime(t)=-\mu_0 x(t)+a_0\,f(x(t-\tau)),\quad \mu_0,a_0\;-\;\text{const},
\end{equation}
which is comprehensively studied in numerous publications (see e.g. \cite{DieSvGSVLWal95,HalSVL93} and further references therein).

Equation (\ref{DDE-Aut}) is used extensively in a broad range  of applications to model various real world phenomena, most notably in 
biology, physiology, and population dynamics (see e.g. \cite{Far17,FriWu92,GlaMac88,LiKua01,Kua93} and further references therein). Some
of such models indeed achieve a better mathematical description and a more accurate prediction of the processes if certain intrinsic 
periodic factors in the environment, such as circadian rhythms and seasonal changes, are taken into account. This naturally leads 
from the initially autonomous DDEs given by Equation~(\ref{DDE-Aut}) to those with periodic coefficients of the form (\ref{DDE-Per}).

In this paper we present some numerical outcomes for DDEs (\ref{DDE-Per}) in the special case of $\mu(t)\equiv0$, when it becomes of the 
form
\begin{equation}\label{DDE-0}
    x^\prime(t)=a(t)f(x(t-\tau)),
\end{equation}
and under the same assumptions of the continuity of  $f$ and $a$ and the $T$-periodicity of the coefficient $a(t)$. 
The corresponding autonomous equation (\ref{DDE-Aut}) has the form  
\begin{equation}\label{DDE_A0}
    x^\prime(t)=a_0f(x(t-\tau)),\quad a_0-\text{const},
\end{equation}
and is even simpler and well studied from the point of view of the existence of periodic solutions and their properties. 
The conditions for the existence of periodic solutions of the latter are 
well known, and can be found in multiple sources, see e.g. \cite{DieSvGSVLWal95,HalSVL93}. To those general conditions mentioned above 
the standard additional ones are that the function $f$ must satisfy the negative feedback condition $x\cdot f(x)<0,\;\forall x\ne0,$ 
be one-sided bounded ($f(x)\le M>0,\; \forall x\in\mathbf{R}$ (or $f(x)\ge -m<0,\; \forall x\in\mathbf{R})$), and the corresponding linearized DDE about the equilibrium $x=0$ must be unstable.

The principal objective we are addressing in this work is to numerically identify and possibly find analytically conditions on the 
coefficient $a(t)$ and the nonlinearity $f$ which would yield the existence of periodic solutions with the same period $T$ as the 
coefficient $a(t)$. A natural approach to start is to consider the existing periodic solutions for the autonomous DDE (\ref{DDE_A0}) 
and follow their changes as the coefficient  $a_0$ is subject to perturbations.  We demonstrate numerically that small periodic
perturbations in the coefficient of the form $a(t)=a_0+\varepsilon a_1(t)$ can lead to induced periodic solutions to DDEs (\ref{DDE-Per})
of the same period. 

Our subsequent and more extensive numerical simulations have lead to precise analytical results and calculations of explicit periodic 
solutions of the same period as that in the coefficient $a(t)$; the corresponding results are contained in the recent papers 
\cite{IvaShe23,IvaShe23-2,IvaShe23-3}.



\section{Preliminaries}
\label{sec:2}

We consider a partial case of the periodic delay differential equation (\ref{DDE-Per}) when $\mu(t)=0$ and $\tau>0$, the coefficient 
$a(t)\ge0$ is non-negative, and the nonlinearity $f$ satisfies the negative feedback assumption $x\cdot f(x)<0, x\ne0.$ 
In addition, the delay $\tau>0$ can always be normalized to the standard one $\tau=1$. Therefore, the equation (\ref{DDE-0}) becomes
\begin{equation}\label{DDE-main}
    x^\prime(t)=a(t)f(x(t-1)).
\end{equation}
In the autonomous case of DDE (\ref{DDE-main}), $a(t)=a_0>0$ unambiguously ensures the negative feedback condition. In this case, the conditions for the existence of periodic solutions are well known.
They can be found in multiple sources, see e.g. monographs \cite{DieSvGSVLWal95,HalSVL93} and further references therein; here we 
present a brief summary of associated results for the sake of readers' convenience. We distinguish between a general case of the 
nonlinearity $f$, and its partial special case when $f$ is an odd function (symmetric).


\subsection{Symmetric (odd) nonlinearity}\label{Sym}
In this subsection we assume that the nonlinearity $f$ is odd, $f(-x)=-f(x),\; \forall x\ne0,$ which case we shall also call the 
``symmetric $f$''. The existence of periodic solutions for the corresponding autonomous DDE (\ref{DDE_A0}) are well-known. 
They are primarily found as the so-called Kaplan-Yorke solutions first described in paper \cite{KapYor74}. 

The following proposition provides conditions for the existence and uniqueness of the symmetric slowly oscillating periodic solution of 
period $4$ ($\tau=1$). One set of conditions guarantees that the solution is stable with the asymptotic phase; the other implies that the
periodic solution is unstable (and therefore not detected numerically). The proposition is a partial case of the more general statements 
in papers \cite{ChoWal88,KapYor74}.

\begin{prop}\label{prop1} (Existence and stability of symmetric periodic solutions).\newline
    Suppose that $f(x)$ is odd, $f(-x)=-f(x)$, satisfying the negative feedback condition $x\cdot f(x)<0,\; \forall x\ne0$, is strictly 
    decreasing on $\mathbf{R}$ with $|f^\prime(0)|>\pi/2$ and $\lim\sup_{x\to\infty}{|f(x)/x|}=\beta<\pi/2$. Then DDE (\ref{DDE_A0}) has
    at least one slowly oscillating periodic solution satisfying the following symmetry conditions: 
    $x(0)=0,\; x(-t)\equiv -x(t),\; \forall t\in\mathbf{R}.$ 
    If, in addition, the derivative $f^\prime (x)$ is monotone on $\mathbf{R_+}$, then:
    \begin{itemize}
        \item[(i)]\; In case when $f^\prime(x)$ is strictly increasing, delay differential equation (\ref{DDE-main}) 
        has a unique slowly oscillating periodic solution of period $4$ which is asymptotically stable with asymptotic phase; 
        \item[(ii)]\; In case when $f^\prime(x)$ is strictly decreasing on $\mathbf{R_+}$, DDE (\ref{DDE-main}) has a unique unstable 
        slowly oscillating periodic solution with period $4$. 
    \end{itemize}
\end{prop}

\subsection{Non-symmetric (general) nonlinearity}\label{N-Sym}
In this subsection we consider a general (non-symmetric) case of the nonlinearity $f$ where only the negative feedback assumption 
is in place. The existence of periodic solutions for the corresponding autonomous DDE (\ref{DDE_A0}) are also well-known. They go back to
the classical work \cite{Jon62b,Jon64}; they are also comprehensively described in monographs \cite{DieSvGSVLWal95,HalSVL93}.

\begin{prop}\label{prop2} (Existence of slowly oscillating periodic solutions). \newline Suppose that $f$ satisfies the negative feedback assumption, $x\cdot f(x) < 0, x\ne0,$ differentiable at zero with $f^\prime(0)>\pi/2$, and is one-sided bounded, $f(x)\le M>0,\; \forall x\in\mathbf{R}$ (or $f(x)\ge -m<0,\; \forall x\in\mathbf{R})$. Then delay differential equation 
    (\ref{DDE-main}) has a non-trivial slowly oscillating periodic solution $x_0(t)$, such that
    $$
    x_0(t)>0\;\forall t\in(0,t_1),\; x_0(t)<0\; \forall t\in(t_1,t_2),\; x_0(t+T)\equiv x_0(t),\; T=t_1+t_2.\;
    $$
\end{prop}
Note that unlike as in Proposition \ref{prop1} the uniqueness of the periodic solution cannot be derived in this more general case. 
In fact there are typical examples when DDE (\ref{DDE_A0}) has multiple periodic solutions \cite{Nus79b,AngBelIvaShe21}. 

\section{Results}

The numerical solutions for the delay differential equation (\ref{DDE-main}) have been calculated using the algorithm described in 
\cite{AngBelIvaShe21}. The algorithm has demonstrated fourth-order global convergence for non-linear systems of delay differential 
equations with discontinuous solutions. The periods of the solutions have been measured and their stability confirmed based on 
persistence of the distances between the consecutive zeros of the obtained solutions with the error less than 0.1\%.


The numerical results in this section are given for the coefficient $a(t)>0$ of the form $a(t)=1+\varepsilon b_0(t)$, where $b_0(t)$ is a fixed $T$-periodic function, and $\varepsilon\in(0,1)$ is a constant parameter. 

We have considered both symmetric and asymmetric functions $f$. For the symmetric case we have computed numerical solutions for the following two non-linearities:
\begin{equation}
f_1\left(x\right)=-C \tanh\left(c x\right),\; C>0, c>0\;
\label{nlin_sym_1}
\end{equation}
and
\begin{equation}
f_2(x)=-\frac{px}{1+x^{2n}},\; p>0, n\ge1, n\in\mathbb{N}.
\label{nlin_sym_2}
\end{equation}

\begin{figure}
\centering
\includegraphics[width=0.49\textwidth] {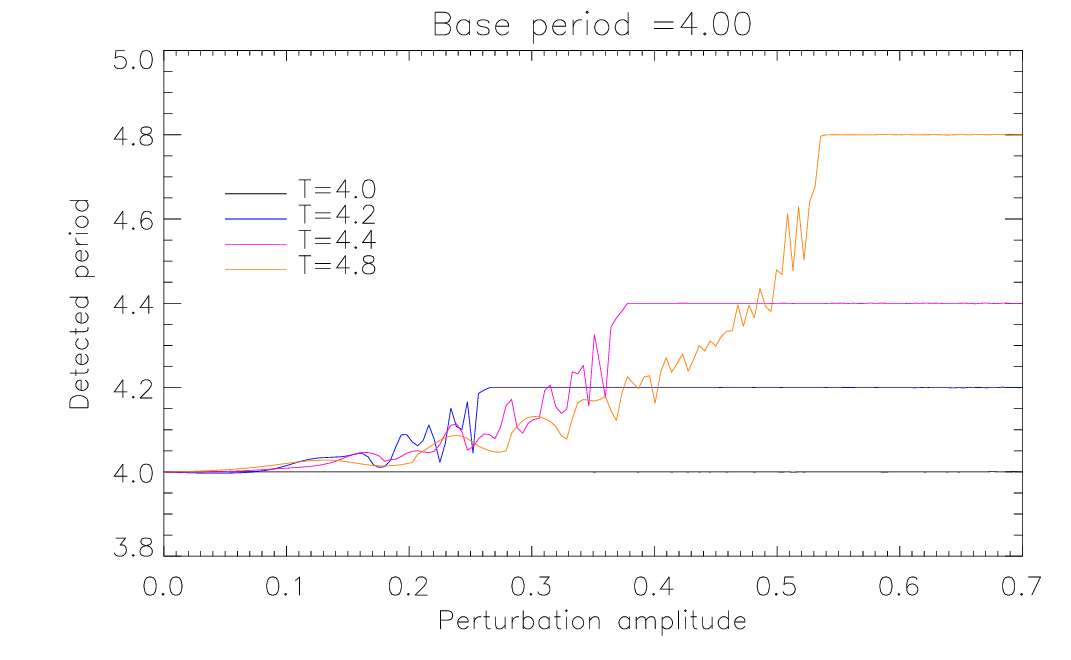}
\includegraphics[width=0.49\textwidth] {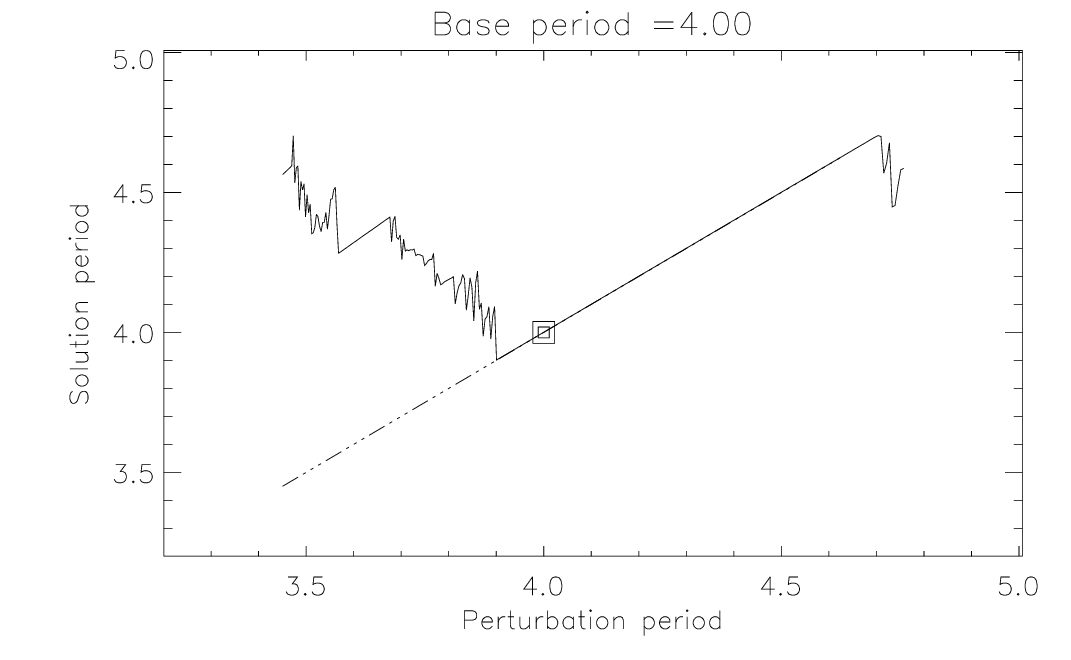}
\caption{Dependence of the detected period in solution of Equation~(\ref{DDE-main}) with the symmetric smooth $f$ given by Equation~(\ref{nlin_sym_1}) on the amplitude of the perturbation $\epsilon$ for three given perturbation periods (left panel) and on the perturbation period $\omega$ for $a=0.5$. The base period of the autonomous equation solution is $T=4$.}
\label{fig:sym_A}
\end{figure}

The numerical data in Figure~\ref{fig:sym_A} is shown for the symmetric case with the first non-linearity $f_1$ when $C=1$ and $c=5$. According to 
Proposition \ref{prop1} for these parametric values DDE (\ref{DDE_A0}) has a unique asymptotically stable slowly oscillating periodic solution with period $T=4$. This solution clearly exists for small perturbation amplitudes. With increasing perturbation amplitude (left panel), the period of the solution grows until it becomes equal to the period of the perturbation.
The effect of varying the perturbation period at a given amplitude of the perturbation is shown in the right panel of Figure~\ref{fig:sym_A}. For a small amplitude $a=0.5$, the range of perturbation periods, for which the period of the solution is equal to the perturbation period, is also narrow in agreement with the left panel of Figure~\ref{fig:sym_A}. Notably, this range is asymmetric. Similar results have been obtained for the non-linearity given by Equation~(\ref{nlin_sym_2}).

\begin{figure}
\centering
\includegraphics[width=0.49\textwidth]{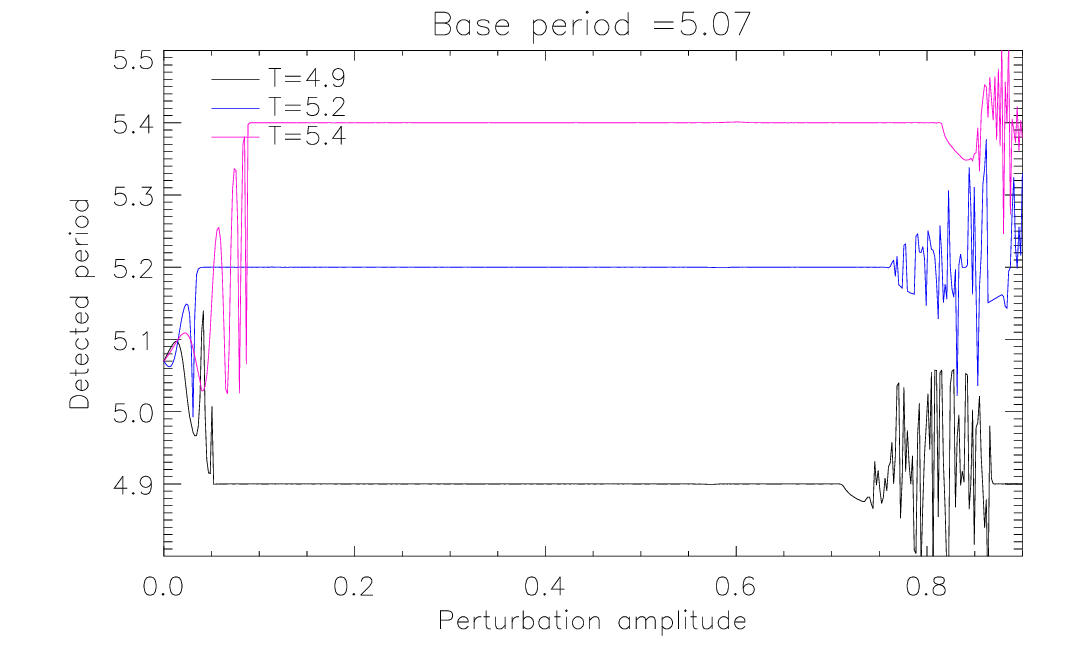}
\includegraphics[width=0.49\textwidth] {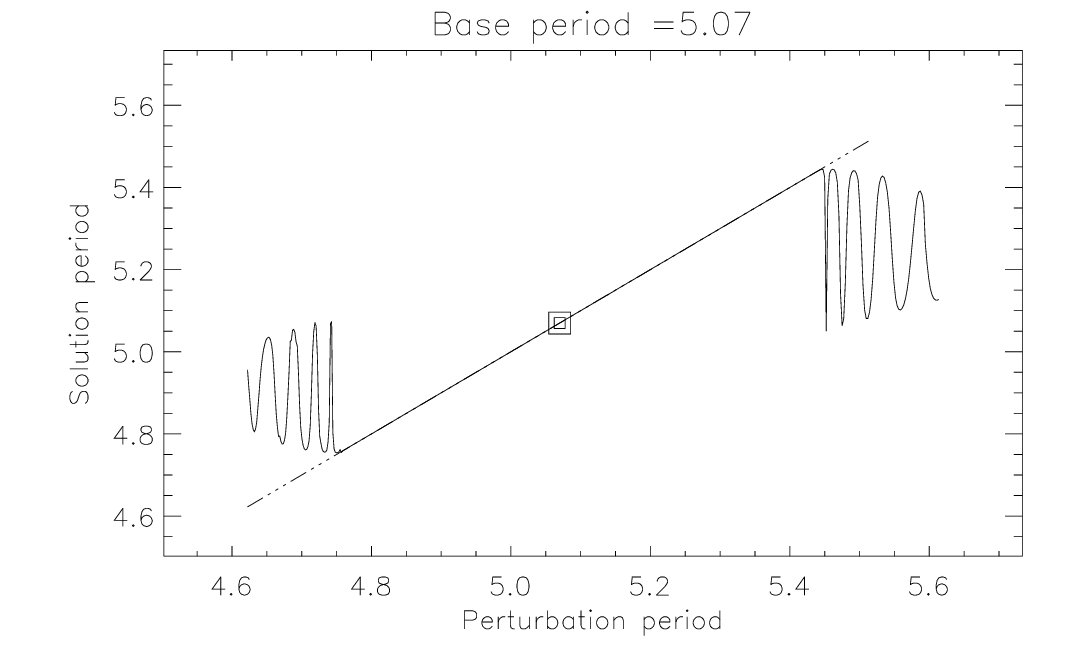}
\caption{Dependence of the detected period in solution of Equation~(\ref{DDE-main}) with $f$ given by Equation~\ref{asym_F_case1} on the amplitude of the perturbation $\epsilon$ for three given perturbation periods (left panel) and on the period $\omega$ of the perturbation with the constant amplitude $a=0.1$ (right panel). The base period of the autonomous equation solution is $T=5.07$.}
\label{fig:asym_A}
\end{figure}

For the asymmetric case, the numerical simulations were run for the following two non-linearities:
\begin{equation}
f_3(x)= A\left[1-\exp\{ax\}\right],\quad A>0,~a>0,
\label{asym_F_case2}
\end{equation}
and
\begin{equation}
    f_4\left(x\right)=   
    \begin{cases}
      -B \arctan\left(b x\right) & \text{if } x<0, \\
      -B \ln\left(1+b x\right) & \text{if } x\geq 0,
    \end{cases} 
\label{asym_F_case1} 
\end{equation}
where $B>0$ and $b>0$.
The latter $f_4$, given by expression ~(\ref{asym_F_case1}), was used to produce the data shown in Figure~\ref{fig:asym_A}. 
The parametric values used are $B=2$ and $b=30$; the period of the corresponding autonomous periodic solution is $T=5.07$. In the first 
case of the non-symmetric $f_3$ (Equation~\ref{asym_F_case2}) analogous computations were performed for $A=2$ and $a=1$  with similar 
outcomes for the those in the case of $f_4$ (they are not depicted by graphs here). The existence of periodic solutions in both 
non-symmetric cases follows from Proposition \ref{prop2}. 

Similar to the symmetric case, as shown in the left panel of Figure~\ref{fig:asym_A}, there is an interval of perturbation amplitudes $\varepsilon$, for which the period of solution is equal to the period of the perturbation. As shown in the right panel of Figure~\ref{fig:asym_A}, there is an interval of perturbation periods for a given perturbation amplitude, for which the period of solution is also equal to the period of perturbation.

\section{Discussion}

In this short paper, we provided a preliminary numerical analysis of a delay differential equation with periodic feedback coefficient. We have demonstrated that stable periodic solutions with the period of the coefficient exist in a wide range of perturbation amplitudes and periods. The next steps will be to extend this study to a broader class of delay-differential equations (e.g. Equation~\ref{DDE-Per}) and to provide an analytical insight into the periodic behaviour of the studied equations.

\begin{acknowledgement}
The authors thank the mathematical research institute MATRIX in Australia where part of this research was performed. Its final version
resulted from collaborative activities of the authors during the workshop "Delay Differential Equations and Their Applications" 
(https://www.matrix-inst.org.au/events/delay-differential-equations-and-their-applications/) held in December 2023. The authors are 
also grateful for the financial support provided for these research activities by Simons Foundation (USA), Flinders University 
(Australia), and the Pennsylvania State University (USA).
A.I.'s research was also supported in part by the Alexander von Humboldt Stiftung (Germany) during his visit to 
Justus-Liebig-Universit\"{a}t, Giessen, in June-August 2023.
\end{acknowledgement}

\end{document}